\documentclass[reqno,draft]{amsart}

\theoremstyle{plain}
\newtheorem{lemma}{Lemma}
\newtheorem{theorem}[lemma]{Theorem}

\newtheorem{definition}[lemma]{Definition}

\theoremstyle{remark}

\newcommand*  {\R} {{\mathbb R}}

\newcommand*{\ang}[1]{\left\langle #1 \right\rangle}

\def\ee{\epsilon}
\def\aa{\alpha}

\def\dd{\delta}

\def\tt{\theta}

\def\la{\lambda}

\def\Gg{\Gamma}
\def\Om{\Omega}

\def\pp{\partial}

\begin{document}

\vskip 0.125in

\title[Salmon's Model]
{Regularity ``in Large"  for the $3D$
 Salmon's Planetary Geostrophic Model of Ocean Dynamics}

\date{December 27, 2010}

\author[C. Cao]{Chongsheng Cao}
\address[C. Cao]
{Department of Mathematics  \\
Florida International University  \\
University Park  \\
Miami, FL 33199, USA.} \email{caoc@fiu.edu}

\author[E.S. Titi]{Edriss S. Titi}
\address[E.S. Titi]
{Department of Mathematics \\
and  Department of Mechanical and  Aerospace Engineering \\
University of California \\
Irvine, CA  92697-3875, USA. {\bf Also:}
 Department of Computer Science and Applied Mathematics \\
Weizmann Institute of Science  \\
Rehovot 76100, Israel.} \email{etiti@math.uci.edu}
\email{edriss.titi@weizmann.ac.il}

\begin{abstract}
It is well known, by now, that the
three-dimensional non-viscous  planetary
geostrophic model, with vertical hydrostatic
balance and horizontal Rayleigh friction, coupled
to the heat diffusion and transport, is
mathematically ill-posed. This is because the
no-normal flow physical boundary condition
implicitly produces an additional boundary
condition for the temperature at the literal
boundary. This additional boundary condition is
different, because of the Coriolis forcing term,
than the no heat flux physical boundary condition.
Consequently, the second order parabolic heat
equation is over determined with two different
boundary conditions. In a previous work we proposed
one remedy to this problem by introducing a
fourth-order artificial hyper-diffusion to the heat
transport equation and proved global regularity for
the proposed model. Another remedy for this problem
was suggested by R.~Salmon by introducing an
additional Rayleigh-like friction term for the
 vertical component of the velocity in the hydrostatic balance
 equation.
In this paper we prove the global, for all time and
all initial data, well-posedness of strong
solutions  to the three-dimensional
 Salmon's planetary geostrophic model of ocean  dynamics. That is, we
 show global existence, uniqueness and
 continuous dependence of the strong solutions
 on initial data for this model.
\end{abstract}

\maketitle



\section{Introduction}   \label{S-1}

The starting point in the derivation of the ocean
circulation models is Boussinesq equations which
are the Navier--Stokes equations with rotation and
a heat transport equation. The global existence of
strong solution to the Navier--Stokes equations,
which is a particular case of the Boussinesq
equations when the temperature is identically zero,
is one of the most challenging problems in applied
analysis. However, geophysicists take advantage of
the shallowness of the oceans and the atmosphere
and introduce  the hydrostatic balance
approximation in the vertical motion. This  in turn
simplifies the Boussinesq model, and leads to the
primitive equations of ocean and atmosphere
dynamics (see, e.g., \cite{LTW92}, \cite{LTW92A},
\cite{Majda03},  \cite{PJ87}, \cite{Richardson},
\cite{S98}, \cite{VA06}
 and references therein).
Further, horizontally, approximations based on the
fast rotation of the earth, and the shallowness of
the atmosphere and ocean imply the smallness of the
Rossby number, which consequently lead to the
geostrophic balance between the Coriolis force and
the horizontal pressure gradient (cf. e.g.,
\cite{GR68}, \cite{PJ87},  \cite{S98}, \cite{VA06}
and references therein). By taking advantage of
these assumptions and other geophysical
considerations several intermediate models
 have been developed and used in numerical studies
of weather
prediction,  long-time climate dynamics and large scale ocean circulation dynamics (see, e.g.,
\cite{BCT88}, \cite{BKK}, \cite{CJ55A}, \cite{CJ55}, \cite{PJ87},
\cite{Richardson}, \cite{SR97}, \cite{SV97A},
\cite{SV97B}, \cite{SD96}, \cite{SH48}, \cite{W88} and references therein).

The planetary geostrophic (PG) model,
 the inviscid and adiabatic form of ``thermocline''  equations,
  of large scale ocean circulation are derived by standard
scaling analysis for gyre--scale oceanic motion
(see \cite{PJ84}, \cite{Ph} , \cite{RS59},
\cite{S98}, \cite{VA06} and \cite{WE59}). They are
given in their simplest dimensionless $\beta-$plane
mid-latitude approximation by the system of
equations:
\begin{eqnarray}
&&\hskip-.8in
 p_x - f  v  = 0, \quad
 p_y + f  u  = 0, \quad  p_z   - T =0, \label{PG-1}  \\
&&\hskip-.8in
u_x +v_y + w_z =0   \label{PG-4} \\
&&\hskip-.8in \pp_t T + u T_x + vT_y + w T_z = \kappa_v T_{zz}\, ,
\label{PG-5}
\end{eqnarray}
in the domain $\Omega = \{ (x,y,z) : (x,y) \in M
\subset \R^2,~ {\mbox{and}}~z \in (-h, 0) \}$. Here
$(u,v,w)$ denotes the velocity field, $p$ is the
pressure, and  $T$ is the temperature, which are
the unknowns.  $ f = f_0 + \beta y$ is the
$\beta-$plane mid-latitude approximation of the
Coriolis force. The first two equations in
(\ref{PG-1})  represent
 the geostrophic balance and the third equation represents
 the hydrostatic balance.
    The diffusive  term, $\kappa_v T_{zz}$
  is  a leading order approximation to the
effect of macro-scale turbulent mixing.
 Based on physical ground
Samelson and Vallis \cite{SV97A} have argued that in closed ocean
basin, with the no-normal-flow boundary conditions, this model can
be solved only in restricted domains which are bounded away from the
lateral boundary, $\partial M \times (-h,0)$. Thus, it cannot be
utilized in the study of the large-scale circulation. Furthermore,
 it has been pointed out numerically in
   \cite{CV86} that arbitrarily  small  linear
disturbances (disturbances that are supported at small spatial
scales) will grow arbitrarily fast when the flow becomes
baroclinically unstable. This nonphysical growth  at small scales is
a signature of mathematical    ill-posedness of this model near
unstable baroclinic mode. Therefore, Samelson and Vallis proposed in
\cite{SV97A} various dissipative schemes to overcome these physical
and numerical
 difficulties. In particular, they propose to add either a linear Rayleigh-like
 drag/friction/damping or a conventional eddy viscosity to the
 horizontal components of the momentum equations, and a
 horizontal diffusion in the thermodynamic equation (subject
 to no-heat-flux at the lateral boundary.)
   The planetary geostrophic (PG) model
with conventional eddy viscosity has been studied
mathematically in \cite{CT03}, \cite{STW98},
\cite{STW00}. In  \cite{CT03}  we show the global
existence and uniqueness of weak and strong
solutions to this 3D viscous PG model. We also
provide rigorous estimates, depending on the
various physical parameters, for the dimension of
its global attractor. In the case where  the
dissipative scheme for the horizontal momentum is
the linear drag Rayleigh friction it is observed
that the second order parabolic PDE that governs
the temperature (the thermodynamic equation) has,
due to the Coriolis force, too many boundary
conditions to be satisfied, and hence it is over
determined and is ill-posed (see, e.g., the
discussion in section \ref{S-2} below,
\cite{CTZ04}, \cite{SV97A}
 and the references therein).  To remedy this
situation  it is argued in \cite{SV97A} that one
would have to add to the thermodynamic equation a
higher order (biharmonic) horizontal diffusion in
order to be able to satisfy both physical boundary
- the no-normal-flow and no-heat-flux boundary
conditions - (cf. e.g., \cite{CTZ04}, \cite{SV97A},
\cite{SV97B}). In \cite{CTZ04}, we introduce,
instead, a new PG model with an appropriate
artificial horizontal ``hyperdiffusion" term, to
the heat equation, which involves the Coriolis
parameters. Under the two natural physical boundary
conditions we are able to prove in \cite{CTZ04} the
global existence and uniqueness of the strong
solutions. Moreover, we also show the existence of
the finite dimensional global attractor.

To overcome the above mentioned non-physical baroclinical instabilities and numerical ill-posedness
Salmon introduced in \cite{S98} the following alternative planetary geostrophic model
 in the cylindrical domain $\Om$:
\begin{eqnarray}
&&\hskip-.8in
\ee \, u - f \, v + p_x   = 0,    \label{SA-1}  \\
&&\hskip-.8in
 \ee\, v + f \, u  + p_y  = 0,   \label{SA-2}   \\
&&\hskip-.8in
\dd \, w + p_z    = T, \label{SA-3}  \\
&&\hskip-.8in
u_x +v_y + w_z =0   \label{SA-4} \\
&&\hskip-.8in \pp_t T -\kappa_h\left(  T_{xx} +  T_{yy}\right)
 - \kappa_v T_{zz} + u T_x + vT_y + w T_z = Q \,, \label{SA-5}
\end{eqnarray}
where $\ee$ and $\dd$ are positive constants representing the linear (Rayleigh friction)
damping coefficients,  and $\kappa_h$ is positive constant
which stand for the horizontal heat diffusivity.  We partition the boundary of $\Om$ into:
\begin{eqnarray}
&&\hskip-.8in
\Gg_u = \{ (x,y,z)  \in \overline{\Om} : z=0 \},  \\
&&\hskip-.8in
\Gg_b = \{ (x,y,z)  \in \overline{\Om} : z=-h \},  \\
&&\hskip-.8in \Gg_s = \{ (x,y,z)  \in \overline{\Om} : (x,y) \in
\partial M, \; -h \leq z \leq 0  \}.
\end{eqnarray}
System (\ref{SA-1})--(\ref{SA-5}) is equipped with the following
boundary conditions -- with no-normal flow and non-heat flux on the
side walls and the bottom (see, e.g., {\bf \cite{LTW92}, \cite{LTW92A},
\cite{PJ87}, \cite{S98}, \cite{SR97}, \cite{SV97A},\cite{SV97B},
\cite{SD96}}):
\begin{eqnarray}
&&\hskip-.8in \mbox{on } \Gg_u:  w=0, \;
\frac{\pp T }{\pp z} +\aa T = 0;
\label{B-1} \\
&&\hskip-.8in \mbox{on }   \Gg_b:  w=0, \; \frac{\pp T }{\pp z}= 0;
\label{B-2} \\
&&\hskip-.8in \mbox{on } \Gg_s: (u, v)  \cdot \vec{n} =0,
\; \frac{\pp T }{\pp \vec{n}} =0,   \label{B-3}
\end{eqnarray}
where  $\vec{n}$ is the normal vector to $\Gg_s$.  In addition, we
supply the system with the initial condition:
\begin{eqnarray}
&&\hskip-.8in T(x,y,z,0) = T_0 (x,y,z).  \label{INIT}
\end{eqnarray}

In this paper we focus on  the question of, and prove, the global regularity
and well-posedness of the $3D$ Salmon's PG model
 (\ref{SA-1})--(\ref{SA-5}) for all time and all initial
 data. We remark that a general discussing concerning the nonlinear system
(\ref{SA-1})--(\ref{INIT}) was presented in
\cite{TTR03}, but without providing any evidence of
its global regularity, a problem that we provide a
positive answer for it in this contribution.

The paper is organized as follows. In section \ref{S-2}, we introduce
our notations and recall some well-known relevant inequalities.  In section
\ref{S-3}  we  show the short-time existence of strong solutions
of system (\ref{SA-1})--(\ref{SA-5}) employing a Galerkin approximation procedure.
Section \ref{S-4} is the main section in which we establish the required
 estimates for proving the global existence and uniqueness of the strong solutions,
 and also show their continuous dependence on the initial data.


\section{Preliminaries}    \label{S-2}

\vskip0.1in

Let us denote by $L^r (\Om)$ and $W^{m, r}(\Om), H^r(\Om)$ the usual $L^r-$Lebesgue and
Sobolev spaces, respectively (cf., \cite{AR75}). We denote by
\begin{equation}
\| \phi\|_r = \left(  \int_{\Om}  |\phi (x,y,z)|^r \; dxdydz \right)^{\frac{1}{r}},
\qquad  \mbox{ for every $\phi \in L^r(\Om)$}.
 \label{L2--1}
\end{equation}
We set
\begin{eqnarray*}
\widetilde{\mathcal{V}} &=&    \left\{ { T \in C^{\infty}(\overline{\Om}):
\left. \frac{\pp T}{\pp z } \right|_{z=-h}= 0; \left. {\left(
\frac{\pp T}{\pp z} + \aa T \right) } \right|_{z= 0}= 0; \; \left.
\frac{\pp T}{\pp \vec{n}} \right|_{\Gg_s}= 0 } \right\},
\end{eqnarray*}
and denote by $V$  the closure spaces of $\widetilde{\mathcal{V}}$ in
$H^1(\Om)$ under the $H^1-$topology. For convenience, we also
introduce the following equivalent norm on $V$:
\begin{equation}
\| \phi\|_V^2 =\kappa_h \|\pp_x \phi (x,y,z)\|_2^2+\kappa_h \|\pp_y \phi (x,y,z)\|_2^2
+ \kappa_v \left( \|\pp_z \phi (x,y,z)\|_2^2+ \aa \|\phi (z=0)\|_{L^2(M)}^2 \right).
 \label{VVV}
\end{equation}
The equivalence of this norm on $V$ to the
$H^1-$norm can be justified thanks to  the
Poincr\'e inequality (\ref{P-2}), below.

 Next, we recall the following
three-dimensional Sobolev and Ladyzhenskaya
inequalities (see, e.g., \cite{AR75}, \cite{CF88},
\cite{GA94}, \cite{LADY})
\begin{eqnarray}
&&\hskip-.68in \| \psi \|_{L^3(\Om)} \leq C_0 \| \psi
\|_{L^2(\Om)}^{1/2} \| \psi
\|_{H^1(\Om)}^{1/2},  \label{SI1}  \\
&&\hskip-.68in \| \psi \|_{L^4(\Om)} \leq C_0 \| \psi
\|_{L^2(\Om)}^{1/4} \| \psi
\|_{H^1(\Om)}^{3/4},
\label{SI2}  \\
&&\hskip-.68in \| \psi \|_{L^6(\Om)} \leq C_0 \| \psi \|_{H^1(\Om)},
\label{SI3}
\end{eqnarray}
for every $\psi\in H^1(\Om).$  Here $C_0$ is a dimensionless positive constant
which might depend on the shape of $M$ and $\Om$  but not on their
sizes. We also introduce the following version of Poincar\'{e} inequality
\begin{eqnarray}
&&\hskip-.68in \| \psi \|_{L^2(\Om)}^2 \leq 2 h \| \psi (z=0) \|_{L^2(M)}^2
+ h^2 \| \psi_z \|_{L^2(\Om)}^2,   \label{P-2}\\
&&\hskip-.68in \| \psi \|_{L^6(\Om)}^6 \leq 2 h \| \psi (z=0) \|_{L^6(M)}^6
+ h^2 \|\psi^2 \; \psi_z \|_{L^2(\Om)}^2.   \label{P-6}
\end{eqnarray}

\vskip0.1in

By solving the linear system (\ref{SA-1})--(\ref{SA-3}) we obtain
\begin{eqnarray}
&&\hskip-.8in
u = - \frac{\ee p_x+fp_y}{\ee^2+f^2},    \label{PG1}  \\
&&\hskip-.8in
v =  \frac{f p_x - \ee p_y}{\ee^2+f^2},     \label{PG2}   \\
&&\hskip-.8in
w =  \frac{T-p_z}{\dd}.    \label{PG3}
\end{eqnarray}
Thanks to (\ref{SA-4}) we have
\begin{eqnarray}
&&\hskip-.8in
- \left[ \left( \frac{\ee p_x+fp_y}{\ee^2+f^2}\right)_x +
\left( \frac{-f p_x + \ee p_y}{\ee^2+f^2} \right)_y+
\left( \frac{p_z-T}{\dd} \right)_z \right]
  = 0.   \label{L-P}
\end{eqnarray}
Using the boundary conditions (\ref{B-1}) and (\ref{B-2}) we infer the
following boundary conditions:
\begin{eqnarray}
&&\hskip-.8in \mbox{on } \Gg_u \;\; \mbox{and} \;\;  \Gg_b:  p_z -
T=0, \quad  \mbox{    and  on } \Gg_s:  \frac{\pp p }{\pp \vec{e}} =0,
\label{PB}
\end{eqnarray}
where $\vec{e} = \frac{\ee \vec{n} + f \vec{k}\times \vec{n}
}{\sqrt{\ee^2+f^2}},$ and $\vec{k}$ is the unit vector of vertical direction.
Notice that by following the techniques developed in
{\bf \cite{Hu-Temam-Ziane}} and {\bf \cite{Ziane}}
(for the case of smooth domains, see,  for example,
 {\bf \cite{LADY}} p. 89, and {\bf \cite{US75}}),
the  three-dimensional second order elliptic
boundary--value problem  (\ref{L-P})--(\ref{PB})
has a unique solution for every given $T$;
moreover, this solution enjoys the following
regularity properties. Taking the $L^2(\Om)$ inner
product of equation (\ref{L-P}) with $p$,
integrating by parts and applying the boundary
conditions (\ref{PB}) and using the Cauchy--Schwarz
inequality, we obtain
\begin{eqnarray}
&&\hskip-.8in
\int_{\Om}\left[ \frac{\ee}{\ee^2+f^2} \left( p_x^2+p_y^2\right) + \frac{p_z^2}{\dd} \right] \; dxdydz
= \frac{1}{\dd} \int_{\Om} T p_z \; dxdydz \leq  \frac{1}{\dd} \|T\|_2 \|p_z\|_2.   \label{H-1}
\end{eqnarray}
Denote by
\begin{eqnarray}
&&\hskip-.8in
0 < F_0 = \min{f} < F_1 = \max{f}.  \label{F-ZERO}
\end{eqnarray}
We observe that the assumption $F_0 > 0$ indicates that we are dealing with a mid-latitude case and away from the equator.
By using (\ref{F-ZERO}) and applying Young's inequality to (\ref{H-1}), we reach
\begin{eqnarray}
&&\hskip-.8in
\int_{\Om}\left[ \frac{\ee}{\ee^2+F_1^2} \left( p_x^2+p_y^2\right) + \frac{p_z^2}{2\dd} \right] \; dxdydz
\leq
\int_{\Om}\left[ \frac{\ee}{\ee^2+f^2} \left( p_x^2+p_y^2\right) + \frac{p_z^2}{2\dd} \right] \; dxdydz
 \leq \frac{1}{2\dd} \|T\|_2^2.   \label{W-L2}
\end{eqnarray}
Furthermore, by (\ref{L-P}) and the above estimate, we have
\begin{eqnarray}
&&\hskip-.8in
\left\|\frac{\ee }{\ee^2+ f^2} \left( p_{xx} + p_{yy}\right) + \frac{p_{zz}}{\dd} \right\|_2 =
 \left\|\frac{\beta p_x (\ee^2-f^2) +2\ee \beta f p_y}{(\ee^2+f^2)^2} + \frac{T_z}{\dd} \right\|_2  \nonumber  \\
&&\hskip-.8in  \leq C
 \left(\frac{\beta (\|p_x\|_2+\|p_y\|_2)}{\ee^2+F_0^2} + \left\|\frac{T_z}{\dd} \right\|_2\right)  \nonumber  \\
&&\hskip-.8in  \leq C
 \left( \frac{\beta (\ee+F_1)}{\ee^{1/2} \dd^{1/2} (\ee^2+F_0^2)} \|T\|_2+ \frac{\|T_z\|_2}{\dd} \right).   \label{H-2}
\end{eqnarray}
As a result of the above and (\ref{PG1})-(\ref{PG3}), we obtain
\begin{eqnarray}
&&\hskip-.8in
\|\ee u\|_2 + \|\ee v\|_2 +
\|\dd \, w\|_2 \leq C (\|\nabla p\|_2 + \|T\|_2 )\leq C \|T\|_2,
\label{W-E-1}
\end{eqnarray}
and
\begin{eqnarray}
&&\hskip-.8in
\|\ee u\|_{H^{1}(\Om)} + \|\ee v\|_{H^{1}(\Om)} +
\|\dd \, w\|_{H^{1}(\Om)} \leq C (\|\nabla p\|_{H^{1}(\Om)}+ \|T\|_{H^{1}(\Om)})
\leq C \|T\|_{H^{1}(\Om)}.
\label{W-E}
\end{eqnarray}

\begin{definition} \label{D-1}
\thinspace Let  $ T_0\in V$, and let $\mathcal{T}$ be a fixed
positive time. \thinspace $(u,v,w,p,T)$  is called a strong solution of {\em
(\ref{SA-1})--(\ref{SA-5})} on the time interval $[0,\mathcal{T}]$ if
\begin{itemize}
\item[1)]
\begin{eqnarray*}
&&  T  \in C([0,\mathcal{T}], V)  \cap L^2 ([0,\mathcal{T}],
H^2(\Om)),     \\
&& T_t  \; \in L^1([0,\mathcal{T}], L^2(\Om)), \\
&& T_t (z=0)  \; \in L^1([0,\mathcal{T}], H^{-1/2}(M)).
\end{eqnarray*}

\item[2)]
$(u, v,w,p)$  satisfies  {\rm (\ref{PG1})--(\ref{PB})}.

\item[3)] Moreover,
{\em (\ref{SA-5})} is satisfied  in the weak sense, namely,
for every $t_0 \in [0, \mathcal{T}]$
\begin{eqnarray}
&&\hskip-0.48in \int_{\Om} T(t) \psi \; dxdydz - \int_{\Om}
T(t_0)\psi   \; dxdydz     \nonumber  \\
&&\hskip-0.48in
+ \int_{t_0}^t \left[ \int_{\Om} \left( \kappa_h   T_{x} \psi_x + \kappa_h T_{y} \psi_y
+\kappa_v T_z \psi_z \right) \; dxdydz  +\kappa_v\, \aa \int_{M} T(z=0) \psi(z=0) \; dxdy \right]\; ds
\label{WEAK}  \\
&&\hskip-0.38in +\int_{t_0}^t \int_{\Om}  \left[ v
 \cdot \nabla T(s) +
w T_z(s)  \right] \psi \; dxdydz
\; ds =\int_{t_0}^t \int_{\Om} Q \psi\;dxdydz \; ds, \nonumber
\end{eqnarray}
for every $\psi \in V,$ and $ t \in
[t_0,\mathcal{T}].$
\end{itemize}
\end{definition}


\section{Short-time Existence of the Strong Solutions} \label{S-3}

In this section we will show the short-time existence of the strong solution
of system (\ref{SA-1})--(\ref{SA-5}).

\begin{theorem} \label{T-WEAK}
Let $Q \in L^2(\Om)$ and   $ T_0\in V$  be given.
Then there exists  a strong solution $(u, v, w, p,
T)$ of system {\em (\ref{SA-1})--(\ref{SA-5})} on
the interval  $[0,\mathcal{T}^{***}]$, where
$\mathcal{T}^{***}$ is a positive time given in
(\ref{T-STAR}), below. Furthermore, $\pp_t T \in
L^2 ([0,\mathcal{T}^{***}]; L^2(\Om))$  and $\pp_t
T (z=0) \in L^2 ([0,\mathcal{T}^{***}];
H^{-1/2}(M))$; and equation (\ref{SA-5}) holds as a
functional  equation in $L^2
([0,\mathcal{T}^{***}]; L^2(\Om))$.
\end{theorem}

\vskip0.05in

\begin{proof}

We will use a Galerkin like procedure to show the existence of the strong
solution for system (\ref{SA-1})--(\ref{SA-5}). First, we will show the existence of the weak solutions.
 Let $\{\phi_k \in V\cap H^2(\Om) \}_{k=1}^{\infty}$ and  $\{\la_k \in \mathbb{R}^+ \}_{k=1}^{\infty}$ be the
eigenfunctions and their corresponding eigenvalues
of the second order elliptic operators
$-\kappa_h\left( T_{xx} + T_{yy}\right)
 - \kappa_v T_{zz}$, subject to the  boundary conditions (\ref{B-1})--(\ref{B-3})
(see, e.g., \cite{LADY}). The eigenvalues are
ordered such that $0 < \la_1 \leq \la_2 \leq
\cdots$; moreover, $\{\phi_k\}_{k=1}^{\infty}$ is
an orthogonal  basis of $L^2(\Om)$. Let $m \in
\mathbb{Z}^+$ be fixed and  $H_m$ be the linear
space generated by
 $\{\phi_k \}_{k=1}^{m}.$ We will denote by $P_m: L^2\rightarrow H_m$, the orthogonal projection
 in $L^2$.
The Galerkin approximating system of order
$m$ that we use for (\ref{SA-1})--(\ref{SA-5}) reads:
\begin{eqnarray}
&&\hskip-.8in
\ee \, u_m - f \, v_m + \pp_x p_m   = 0,    \label{GSA-1}  \\
&&\hskip-.8in
 \ee\, v_m + f \, u_m  + \pp_y p_m  = 0,   \label{GSA-2}   \\
&&\hskip-.8in
\dd \, w_m + \pp_z p_m    = T_m,   \label{GSA-3}  \\
&&\hskip-.8in
\pp_x u_m + \pp_y v_m +\pp_z w_m =0   \label{GSA-4} \\
&&\hskip-.8in \pp_t T_m -\kappa_h \left(  \pp_{xx} T_m  +  \pp_{yy} T_m \right)
 - \kappa_v \pp_{zz} T_m  + P_m \left[ u_m \pp_x T_m + v_m \pp_y T_m + w_m \pp_z T_m \right] = P_m Q,
\, \label{GSA-5}
  \\
&&\hskip-.8in
T_m(x,y,z,0) = P_m T_0 (x,y,z),
\label{GSA-6}
\end{eqnarray}
where $T_m =  \sum_{k=1}^m a_k(t) \phi_k(x,y,z)$,
and $(u_m, v_m,  w_m, p_m)$ is the solution of the
system (\ref{GSA-1})--(\ref{GSA-4}) under boundary
condition $ \left. w_m\right|_{z=0}=\left.
w_m\right|_{z=-h}=0; \left. (u_m, v_m)\cdot \vec{n}
\right|_{\Gg_s}    =0.$  Based on discussion in the
previous section, equation (\ref{GSA-5}) is an ODE
system with the unknown $a_k(t), k=1,\cdots, m.$
Furthermore, it is easy to check that the vector
field in equation (\ref{GSA-5}) is locally
Lipschitz with respect to $a_k(t), k=1,\cdots, m,$
since it is quadratic. Therefore, there is a unique
solution $a_k(t), k=1,\cdots, m,$ to equation
(\ref{GSA-5}) for a short interval of time $[0,
\mathcal{T}^*_m]$. Let $[0, \mathcal{T}^{**}_m)$ be
the maximal interval of existence for system
(\ref{GSA-1})--(\ref{GSA-6}). We will focus our
discussion below on the interval $[0,
\mathcal{T}^{**}_m)$, and will show that
$\mathcal{T}^{**}_m= +\infty.$

By taking the $L^2(\Om)$ inner product of equation
(\ref{GSA-5}) with $T_m$, we obtain
\begin{eqnarray}
&&\hskip-.28in \frac{1}{2} \frac{ d \|T_m\|_2^2 }{d t}+ \kappa_h
\left( \|\pp_x T_m\|_2^2 +\|\pp_y T_m\|_2^2\right)
+ \kappa_v \left(\|\pp_z T_m\|_2^2 +\aa \|T_m(z=0)\|_2^2  \right)  \\
&&\hskip-.228in
 + \int_{\Om} \left[ u_m \pp_x T_m + v_m \pp_y T_m + w_m
\pp_z T_m \right]  \, T_m  \, dxdydz =\int_{\Om} Q \, T_m \, dxdydz.
\label{L2}
\end{eqnarray}
It is easy to show by integrating by parts and
by using the relevant boundary conditions
(\ref{B-1})--(\ref{B-3}) that
\begin{equation}
\int_{\Om} \left[ u_m \pp_x T_m + v_m \pp_y T_m +
w_m \pp_z T_m \right]  \, T_m  \, dxdydz  =0.
\label{EST-1}
\end{equation}
Furthermore, by the Cauchy--Schwarz inequality and (\ref{P-2}) we have
\begin{eqnarray*}
&&\hskip-.18in
\left| { \int_{\Om} Q\, T_m \;  dxdydz } \right|
\leq  \| Q \|_{2} \| T_m \|_{2} \\
&&\hskip-.18in
\leq \frac{1}{\sqrt{\la_1}} \| Q \|_{2} \left[  \kappa_h
\left( \|\pp_x T_m\|_2^2 +\|\pp_y T_m\|_2^2\right)
+ \kappa_v \left(\|\pp_z T_m\|_2^2 +\aa \|T_m(z=0)\|_2^2  \right)  \right],
\end{eqnarray*}
where $\la_1$ is the first eigenvalue discussed above.
From the above estimates, we obtain
\begin{eqnarray}
&&\hskip-.18in
\frac{ d \|T_m\|_2^2 }{d t}
+ \kappa_h
\left( \|\pp_x T_m\|_2^2 +\|\pp_y T_m\|_2^2\right)
+ \kappa_v \left(\|\pp_z T_m\|_2^2 +\aa \|T_m(z=0)\|_2^2  \right)
\leq \frac{\|Q\|_2^2}{\la_1}.     \label{T-STAR1}
\end{eqnarray}
Consequently, we have,
\begin{eqnarray*}
&&\hskip-.18in
\frac{ d \|T_m\|_2^2 }{d t}
+ \la_1 \|T_m\|_2^2
\leq \frac{\|Q\|_2^2}{\la_1}.
\end{eqnarray*}
Thanks to Gronwall inequality, we conclude that
\begin{eqnarray}
&&\hskip-.3in \|T_m(t)\|_2^2 \leq \|T_0\|_2^2 \; e^{-\la_1 \; t} + \frac{\|Q\|_2^2}{\la_1^2},
 \label{EST-L-2}
\end{eqnarray}
for every $t \in [0, \mathcal{T}^{**}_m)$. From the
above, we conclude that $T_m(t)$ must exist
globally, i.e., $\mathcal{T}^{**}_m = +\infty.$
Therefore, for any given $\mathcal{T} >0$ and any
$t \in [0,\mathcal{T}]$, we have
\begin{eqnarray}
&&\hskip-.3in \|T_m(t)\|_2^2 \leq \|T_0\|_2^2 \; e^{-\la_1 \; t} +
\frac{\|Q\|_2^2}{\la_1^2}. \label{WL-2}
\end{eqnarray}
Furthermore, by integrating (\ref{T-STAR1}) with
respect the time variable over the interval
$[0,t]$,  for  $t \in [0,\mathcal{T}]$, and by
(\ref{WL-2}), we get
\begin{eqnarray}
&&\hskip-.3in \int_{0}^t
\left[ \kappa_h
\left( \|\pp_x T_m (s)\|_2^2 +\|\pp_y T_m(s)\|_2^2\right)
+ \kappa_v \left(\|\pp_z T_m(s)\|_2^2 +\aa \|T_m(z=0)(s)\|_2^2  \right)
\right]
 \; ds  \nonumber   \\
 &&\hskip-.3in \leq \|T_0\|_2^2+ \frac{\|Q\|_2^2 \; t}{\la_1^2}.
\label{WLL-2}
\end{eqnarray}
As a result of all the above we have established
that $T_m$ exists globally in time, and that it is
uniformly bounded, with respect to $m$, in the
$L^{\infty} ([0,\mathcal{T}]; L^2(\Om))$ and $L^2
([0,\mathcal{T}]; V)$ norms.

Next, and similar to the theory of $3D$ Navier--Stokes
equations (see, e.g., \cite{CF88} and \cite{TT84}),
let us show that $\pp_t T_m$ is uniformly bounded, with respect to
 $m$, in the  $L^{\frac{4}{3}} ([0,\mathcal{T}]; V^{\prime})$ norm,
 where $V^{\prime}$ is the dual space of $V$.   {F}rom
(\ref{GSA-5}), we have, for every $\psi \in V$
\[
\ang{\pp_t T_m, \psi} = \ang{ P_m Q  +\kappa_h (\pp_{xx} T_m
+\pp_{yy} T_m) + \kappa_v \pp_{zz} T_m - P_m\left[u_m \pp_x T_m +
v_m \pp_y T_m+ w_m \pp_z T_m \right], \psi}.
\]
Here, $\ang{\cdot, \cdot}$ is the dual action of $V^{\prime}$.  It
is clear that
\begin{eqnarray}
&&\hskip-.5in \left| \ang{P_m Q, \psi}\right| \leq \|Q\|_2
\|\psi\|_2, \label{L2-1}
\end{eqnarray}
and  by integration by parts and using boundary condition (\ref{B-1})--(\ref{B-3}), we have
\begin{eqnarray}
&&\hskip-.5in \left| \ang{  \kappa_h (\pp_{xx} T_m +\pp_{yy} T_m) +
\kappa_v \pp_{zz} T_m,  \psi} \right| \leq C \|T_m \|_V \;
\|\psi\|_V, \label{L2-2}
\end{eqnarray}
recall that $\|\cdot\|_V$ is defined in (\ref{VVV}).
Next, let us get an estimate for
\begin{eqnarray*}
&&\hskip-0.25in \left| \ang{ P_m\left[u_m \pp_x T_m +
v_m \pp_y T_m+ w_m \pp_z T_m  \right], \psi} \right|  \\
&&\hskip-0.18in = \left| \int_{\Om}  \left[u_m \pp_x T_m + v_m \pp_y
T_m+ w_m \pp_z T_m  \right]\; \psi_m \; dxdydz \right|,
\end{eqnarray*}
where $\psi_m = P_m \psi.$ Thus, by integration by parts
and using (\ref{GSA-4}), (\ref{W-E-1}),
(\ref{W-E}) and relevant  boundary conditions, we obtain
\begin{eqnarray}
&&\hskip-0.25in \left| \ang{ P_m\left[u_m \pp_x T_m + v_m \pp_y T_m+
w_m \pp_z T_m \right], \psi} \right|     \nonumber \\
&&\hskip-0.18in = \left| \int_{\Om}  \left[u_m \pp_x \psi_m + v_m
\pp_y \psi_m+ w_m \pp_z \psi_m  \right]\;  T_m \; dxdydz \right|    \nonumber  \\
&&\hskip-0.8in \leq C \left[ \| u_m \|_4 +\| v_m \|_4 +\| w_m \|_4
\right] \; \| T_m \|_4 \;  \| \nabla \psi_m \|_2       \\
&&\hskip-0.8in \leq C \left(  \| u_m \|_2^{1/4} \|u_m \|_{H^1}^{3/4} + \| v_m \|_2^{1/4} \|v_m \|_{H^1}^{3/4} + \| w_m \|_2^{1/4} \|w_m \|_{H^1}^{3/4}\right)  \| T_m \|_2^{1/4} \|T_m \|_{H^1}^{3/4}
\| \nabla \psi_m \|_2  \nonumber  \\
&&\hskip-0.8in \leq C  \left( \| T_m \|_2^2 + \| T_m \|_2^{1/2} \|T_m \|_V^{3/2} \right) \;
\| \nabla \psi \|_2. \label{L2-3}
\end{eqnarray}
Therefore, by the estimates (\ref{L2-1})--(\ref{L2-3}),   we have
\begin{eqnarray*}
&&\hskip-.25in
\left| \ang{ \pp_t T_m, \psi } \right|
\leq C \left(\|Q\|_2 + \|T_m \|_V +\| T_m \|_2^2+ \| T_m \|_2^{1/2} \|T_m \|_V^{3/2} \right) \|  \psi\|_V.
\end{eqnarray*}
Thus, we have
\begin{equation}
\int_0^t \|\pp_t T_m (t)
\|_{V^\prime}^{\frac{4}{3} } dt
\leq C \left( \|Q\|_2\; t^{4/3} +  \|T_0\|_2^2+ \frac{\|Q\|_2^2 \; t}{\la_1^2} \right).
\label{DL2}
\end{equation}
Therefore, $\pp_t T_m$ is uniformly bounded,
 with respect to $m$, in the  $L^{\frac{4 }{3} } ([0,\mathcal{T}]; V^{\prime})$ norm.
Thanks to (\ref{WL-2}), (\ref{WLL-2}) and (\ref{DL2}),
one can apply the Aubin's compactness Theorem
(cf., for example, {\bf \cite{CF88},
\cite{TT84}})  and extract a subsequence $\{ T_{m_j} \}$ of
$\{ T_m \}$
and a subsequence $\{ \pp_t T_{m_j} \}$
of $\{ \pp_t T_m \}$;
which converge to
$T \in L^{\infty} ([0,\mathcal{T}]; L^2(\Om))
\cap L^2 ([0,\mathcal{T}]; V)$
and $\pp_t T \in L^{\frac{4}{3} }
([0,\mathcal{T}]; V^{\prime})$, respectively,
 in the following sense:
\[
\left\{ {
\begin{array}{ll}
\displaystyle{T_{m_j} \rightarrow T}
& \displaystyle{ \mbox{in } L^{2} ([0,\mathcal{T}]; L^2(\Om))}  \quad \mbox{strongly}; \\
\displaystyle{T_{m_j} \rightarrow T}
& \displaystyle{ \mbox{in } L^{\infty} ([0,\mathcal{T}]; L^2(\Om))}  \quad \mbox{weak-star}; \\
\displaystyle{ T_{m_j} \rightarrow T}
&\displaystyle{ \mbox{in } L^2 ([0,\mathcal{T}]; H^1(\Om)) \quad \mbox{weakly};} \\
\displaystyle{ \pp_t T_{m_j} \rightarrow
\pp_t T}
&\displaystyle{\mbox{in } L^{\frac{4}{3} } ([0,\mathcal{T}]; V^{\prime} )
\quad \mbox{weakly}.}
\end{array} } \right.
\]
Moreover, from (\ref{GSA-1})--(\ref{GSA-4}) (see
also (\ref{SA-1})--(\ref{SA-4})) we observe that
 $\{u_m, v_m, w_m \}$  depend linearly on $T_m$. Therefore,
 the elliptic estimates (\ref{W-E-1}) and (\ref{W-E})
 imply, thanks to (\ref{WL-2}) and (\ref{WLL-2}), uniform bounds,
 with respect to $m$, for $\{u_m , v_m , w_m  \}$
 in $L^{\infty} ([0,\mathcal{T}]; L^2(\Om))$ and
 $L^{2} ([0,\mathcal{T}]; H^1(\Om))$, respectively. Therefore,
 we can extract a subsequence of $\{ u_{m_j}, v_{m_j} w_{m_j}
 \}$, corresponding to the readily established
 converging subsequence for the temperature $\{
 T_{m_j}\}$, which will be also labeled $\{ u_{m_j}, v_{m_j} w_{m_j}
 \}$, that converges to $\{
u, v, w\}$ weak-star in $L^{\infty}
([0,\mathcal{T}]; L^2(\Om))$, and weakly in $L^{2}
([0,\mathcal{T}]; H^1(\Om))$. By passing to the
limit, one can show as in the case of
Navier--Stokes equations (see, for example, {\bf
\cite{CF88}, \cite{TT84}}) that $T$ also satisfies
(\ref{WEAK}). In other words, $T$ is a weak
solution of the system (\ref{SA-1})--(\ref{SA-5}).

By taking the $L^2(\Om)$ inner product of equation
(\ref{GSA-5}) with $-\kappa_h (\pp_{xx} T_m +\pp_{yy} T_m) -
\kappa_v \pp_{zz} T_m$, we obtain
\begin{eqnarray*}
&&\hskip-.28in \frac{1}{2} \frac{ d}{dt} \left[ \kappa_h
\left( \|\pp_x T_m\|_2^2 +\|\pp_y T_m\|_2^2\right)
+ \kappa_v \left(\|\pp_z T_m\|_2^2 +\aa \|T_m(z=0)\|_2^2  \right) \right]  +  \|\kappa_h (\pp_{xx} T_m +\pp_{yy} T_m) +
\kappa_v \pp_{zz} T_m \|_2^2 \\
&&\hskip-.28in = \int_{\Om} \left(  Q - u_m \pp_x T_m + v_m \pp_y T_m + w_m
\pp_z T_m \right) \, \left(\kappa_h (\pp_{xx} T_m +\pp_{yy} T_m) +
\kappa_v \pp_{zz} T_m\right) \, dxdydz   \\
&&\hskip-.28in
\leq  \left( \|Q\|_2 \,  +  \| u_m \|_6 \|\pp_x T_m \|_3 + \|v_m\|_6 \|\pp_y T_m \|_3+ \|w_m\|_6
\|\pp_z T_m \|_3 \right) \, \|\kappa_h (\pp_{xx} T_m +\pp_{yy} T_m) +
\kappa_v \pp_{zz} T_m \|_2   \\
&&\hskip-.28in
\leq  \left( \|Q\|_2 \,  +  C \| T_m \|_6 \|\nabla T_m \|_3  \right) \, \|\kappa_h (\pp_{xx} T_m +\pp_{yy} T_m) +
\kappa_v \pp_{zz} T_m \|_2 \\
&&\hskip-.28in
\leq  \left( \|Q\|_2 + C\|T_m\|_2^2\right) \, \|\kappa_h (\pp_{xx} T_m +\pp_{yy} T_m) +
\kappa_v \pp_{zz} T_m \|_2  \\
&&\hskip-.228in + C\left[ \kappa_h
\left( \|\pp_x T_m\|_2^2 +\|\pp_y T_m\|_2^2\right)
+ \kappa_v \left(\|\pp_z T_m\|_2^2 +\aa \|T_m(z=0)\|_2^2  \right) \right]^{\frac{3}{2}}  \, \|\kappa_h (\pp_{xx} T_m +\pp_{yy} T_m) +
\kappa_v \pp_{zz} T_m \|_2^{\frac{3}{2}}.
\end{eqnarray*}
Therefore, applying the Cauchy--Schwarz inequality
and Young's inequality to the above estimate,  we
obtain
\begin{eqnarray}
&&\hskip-.28in \frac{ d}{dt} \left[ \kappa_h
\left( \|\pp_x T_m\|_2^2 +\|\pp_y T_m\|_2^2\right)
+ \kappa_v \left(\|\pp_z T_m\|_2^4 +\aa \|T_m(z=0)\|_2^2  \right) \right]  +  \|\kappa_h (\pp_{xx} T_m +\pp_{yy} T_m) +
\kappa_v \pp_{zz} T_m \|_2^2  \nonumber  \\
&&\hskip-.28in
\leq  \|Q\|_2^2 + C \|T_m\|_2^4 + C\left[ \kappa_h
\left( \|\pp_x T_m\|_2^2 +\|\pp_y T_m\|_2^2\right)
+ \kappa_v \left(\|\pp_z T_m\|_2^2 +\aa \|T_m(z=0)\|_2^2  \right) \right]^{6}.
\label{T-STAR-1}
\end{eqnarray}
Consequently, we have
\begin{eqnarray*}
&&\hskip-.3in
\kappa_h
\left( \|\pp_x T_m\|_2^2 +\|\pp_y T_m\|_2^2\right)
+ \kappa_v \left(\|\pp_z T_m\|_2^2 +\aa \|T_m(z=0)\|_2^2  \right)   \\
&&\hskip-.28in
 \leq \frac{\kappa_h
\left( \|\pp_x T_0\|_2^2 +\|\pp_y T_0\|_2^2\right)
+ \kappa_v \left(\|\pp_z T_0\|_2^2 +\aa \|T_0(z=0)\|_2^2  \right)}{\left(1- C \, t \; \left( \|T_0\|_2^4 +
\|Q\|_2^2 \right)\,  \left[ \kappa_h
\left( \|\pp_x T_0\|_2^2 +\|\pp_y T_0\|_2^2\right)
+ \kappa_v \left(\|\pp_z T_0\|_2^2 +\aa \|T_0(z=0)\|_2^2  \right) \right]\right)^{1/2}}.
\end{eqnarray*}
Therefore, for every $t \in [0, \mathcal{T}^{***}]$, where
\begin{eqnarray}
&&\hskip-.3in
\mathcal{T}^{***} := \frac{1}{4C \left( \left( \|T_0\|_2^4 +
\|Q\|_2^2 \right)\,  \left[ \kappa_h
\left( \|\pp_x T_0\|_2^2 +\|\pp_y T_0\|_2^2\right)
+ \kappa_v \left(\|\pp_z T_0\|_2^2 +\aa \|T_0(z=0)\|_2^2  \right) \right]\right)},
 \label{T-STAR}
\end{eqnarray}
  we have
\begin{eqnarray}
&&\hskip-.3in
\kappa_h
\left( \|\pp_x T_m\|_2^2 +\|\pp_y T_m\|_2^2\right)
+ \kappa_v \left(\|\pp_z T_m\|_2^2 +\aa \|T_m(z=0)\|_2^2  \right)  \nonumber  \\
&&\hskip-.28in
 \leq 2 \left[\kappa_h
\left( \|\pp_x T_0\|_2^2 +\|\pp_y T_0\|_2^2\right)
+ \kappa_v \left(\|\pp_z T_0\|_2^2 +\aa \|T_0(z=0)\|_2^2  \right) \right].
 \label{EST-L-2-m}
\end{eqnarray}
Moreover, by integrating (\ref{T-STAR-1}) we obtain
\begin{eqnarray}
&&\hskip-.3in \int_{0}^t \|\kappa_h (\pp_{xx} T_m(s) +\pp_{yy} T_m (s)) +
\kappa_v \pp_{zz} T_m (s)\|_2^2 \; ds  \nonumber  \\
&&\hskip-.28in
\leq \kappa_h
\left( \|\pp_x T_0\|_2^2 +\|\pp_y T_0\|_2^2\right)
+ \kappa_v \left(\|\pp_z T_0\|_2^2 +\aa \|T_0(z=0)\|_2^2  \right) +
\|Q\|_2^2 \; t + C \left(\|T_0\|_2^2 \; e^{-\la_1 \; t} + \frac{\|Q\|_2^2}{\la_1^2} \right) \; t \nonumber  \\
&&\hskip-.228in
+ C\left[ \kappa_h
\left( \|\pp_x T_0\|_2^2 +\|\pp_y T_0\|_2^2\right)
+ \kappa_v \left(\|\pp_z T_0\|_2^2 +\aa \|T_0(z=0)\|_2^2  \right) \right]^{6} \; t,   \qquad t \in [0, \mathcal{T}^{***}].
\label{WLL-21}
\end{eqnarray}
Notice that $T_m$ exists, globally. What we have
just proved is that the $L^2
([0,\mathcal{T}^{***}]; H^2(\Om))$ norm of $T_m$ is
bounded uniformly with respect to $m$. As a result
of all the above we have  $T_m$ exists, at least,
on $[0, \mathcal{T}^{***}]$ and is uniformly
bounded, with respect to $m$, in $L^{\infty} ([0,
\mathcal{T}^{***}]; V)$ and $L^2
([0,\mathcal{T}^{***}]; H^2(\Om))$ norms.
Furthermore, and as for the theory of the
Navier-Stokes equations (see, for example, {\bf
\cite{CF88}, \cite{TT84}}), we can use the above
bounds (\ref{EST-L-2-m}) and (\ref{WLL-21}) to show
that the $L^{2} ([0, \mathcal{T}^{***}]; L^2(\Om))$
norm of $\pp_t T_m$ and the $L^2
([0,\mathcal{T}^{***}]; H^{-1/2} (M))$ norm of
$\pp_t T_m(z=0)$ are uniformly bounded with respect
to $m$. Passing to the limits, we conclude that
there is a strong solution to system
(\ref{SA-1})--(\ref{SA-5}), at least, on $[0,
\mathcal{T}^{***}]$. Furthermore, this strong
solution enjoys the following properties:
\begin{eqnarray}
\pp_t T \in L^2 ([0,\mathcal{T}^{***}]; L^2(\Om))   \quad  \mbox{and}  \quad
\pp_t T (z=0) \in L^2 ([0,\mathcal{T}^{***}]; H^{-1/2}(M)).   \label{DER}
\end{eqnarray}

The above regularity estimates are sufficient to
complete the proof of Theorem \ref{T-WEAK},
following standard techniques from the theory of
the Navier--Stokes equations (see, e.g.,
\cite{CF88} and \cite{TT84}). Furthermore, as a
consequence of the above estimates, in particular
those implying (\ref{DER}), we conclude that
equation (\ref{SA-5}) holds as a functional
equation in $L^2 ([0,\mathcal{T}^{***}];
L^2(\Om))$.
\end{proof}

\section{Global Existence and Uniqueness of the Strong Solutions} \label{S-4}

In the previous section we have established the
short-time existence of the strong solution
to system (\ref{SA-1})--(\ref{SA-5}).
In this section  we will  show the global existence and
uniqueness, i.e. global regularity, of strong solutions to the
system (\ref{SA-1})--(\ref{SA-5}), and their continuous dependence on initial data.

\begin{theorem} \label{T-MAIN}
Let $Q \in L^2(\Om)$,   $ T_0\in V$ and $\mathcal{T}>0,$ be given.
Then there exists  a unique strong solution $(u,v,w,p, T)$ of the system
{\em (\ref{SA-1})--(\ref{SA-5})}, on the interval $[0,\mathcal{T}]$,
which depends continuously on the initial data in the sense specified in equation (\ref{CCC}) below.
\end{theorem}

\vskip0.05in

\begin{proof}
Denote by  $(u,v,w,p,T)$  the strong solution
corresponding to the initial data $T_0$ with
maximal interval of existence $[0,\mathcal{T}_*)$,
that has been established in Theorem \ref{T-WEAK}.
We will show  that $\mathcal{T}_* = \infty.$ To
show this we assume by contradiction that
$\mathcal{T}_* < \infty$. Consequently, it is clear
that
\[
\limsup_{t \to \mathcal{T}_*^{-}} \|  T (t)
\|_{H^1(\Om)}    = \infty,
\]
because, otherwise,  and by virtue of Theorem
\ref{T-WEAK}, the solution can be extended beyond
the maximal  time of existence, $\mathcal{T}_*$.
Next, we will show that $\|  T (t) \|_{H^1(\Om)}$
is bounded uniformly on the interval
$[0,\mathcal{T}_*)$. In what follows we will focus
our discussion and estimates  on the finite maximal
interval of existence $[0,\mathcal{T}_*)$.

\subsection{$L^2$ estimates}

As a result of Theorem \ref{T-WEAK}, equation
(\ref{SA-5}) holds  in $L_{\mbox
{loc}}^2([0,\mathcal{T}_*);L^2(\Om))$, therefore we
can take the inner product of equation (\ref{SA-5})
with $T$, in $L^2(\Om)$, and obtain
\begin{eqnarray*}
&&\hskip-.68in \frac{1}{2} \frac{d \|T\|_2^2}{dt} + \kappa_h
\left( \|\pp_x T\|_2^2 +\|\pp_y T\|_2^2\right)
+ \kappa_v \left( \|\pp_z T\|_2^2 +\aa \|T(z=0)\|_2^2  \right) \\
&&\hskip-.65in = \int_{\Om} Q T \; dxdydz  -\int_{\Om} \left( u
\pp_x T +v \pp_y T+ w \pp_z T \right) T \; dxdydz.
\end{eqnarray*}
After integrating by parts we get
\begin{eqnarray}
&&\hskip-.065in   \int_{\Om} \left( u \pp_x T +v  \pp_y T +
w  \pp_z T  \right) T  \; dxdydz =0. \label{DT-9}
\end{eqnarray}
As a result of the above we conclude
\begin{eqnarray*}
&&\hskip-.68in \frac{1}{2} \frac{d \|T \|_2^2}{dt} + \kappa_h
\left( \|\pp_x T \|_2^2 +\|\pp_y T \|_2^2\right)
+ \kappa_v \left(  \|\pp_z T \|_2^2 +\aa \|T(z=0)\|_2^2 \right) \\
&&\hskip-.65in =\int_{\Om} Q T  \; dxdydz \leq \|Q\|_2 \;
\|T \|_2.
\end{eqnarray*}
Using (\ref{P-2}) and  the Cauchy--Schwarz inequality we obtain
\begin{eqnarray}
&&\hskip-.68in  \frac{d \|T \|_2^2}{dt} + 2 \kappa_h
\left( \|\pp_x T \|_2^2 +\|\pp_y T \|_2^2\right)
+ \kappa_v \left( \|\pp_z T \|_2^2 +\aa \|T(z=0)\|_2^2 \right)\\
&&\hskip-.65in  \leq  \left( \frac{h^2}{\kappa_v}  + \frac{2h}{\aa}\right) \|Q\|^2_2.
\label{T_E}
\end{eqnarray}
By the inequality (\ref{P-2}) and thanks to Gronwall inequality the
above gives
\begin{eqnarray}
&&\hskip-.68in  \|T \|_2^2  \leq
  e^{-\; \frac{t}{2(h^2 + h/\aa)}} \|T_0\|_2^2 +
(2 h^2 + 2h/\aa)^2 \|Q\|^2_2,   \label{T-2}
\end{eqnarray}
for are $t\in[0,\mathcal{T}_*)$. Moreover, we also
have
\begin{eqnarray}
&&\hskip-.68in   \int_0^t \left[   \kappa_h \left( \|\pp_x T \|_2^2
+\|\pp_y T \|_2^2\right) + \kappa_v \left( \|\pp_z T \|_2^2 +\aa \|T(z=0)\|_2^2 \right) \right]\; ds
\nonumber \\
&&\hskip-.65in \leq  2(h^2 + \frac{h}{\aa}) \|Q\|^2_2 \; t +
 e^{-\; \frac{t}{2(h^2 + h/\aa)}} \|T_0\|_2^2 +
(2 h^2 + 2h/\aa)^2 \|Q\|^2_2,    \label{T-2I}
\end{eqnarray}
for are $t\in[0,\mathcal{T}_*)$.

We remark that estimates (\ref{T-2}) and
(\ref{T-2I}) also follow directly from (\ref{WL-2})
and (\ref{WLL-2}), respectively.


\subsection{$L^6$ estimates}
Recall from Theorem \ref{T-WEAK} that $T \in
L_{\mbox {loc}}^{\infty}([0, \mathcal{T}_*),
H^1(\Om)) \cap L_{\mbox {loc}}^2([0,
\mathcal{T}_*), H^2(\Om))$, therefore $|T|^4 T \in
L_{\mbox {loc}}^2([0,\mathcal{T}_*);L^2(\Om))$.
Since by Theorem \ref{T-WEAK} equation (\ref{SA-5})
holds in $L_{\mbox
{loc}}^2([0,\mathcal{T}_*);L^2(\Om))$ we can take
the  inner product of the equation (\ref{SA-5}), in
$L^2(\Om)$, with  $|T|^4 T$ to get
\begin{eqnarray*}
&&\hskip-.168in \frac{1}{6} \frac{d \| T \|_{6}^{6} }{d t} +
 5 \int_{\Om} \left[ \kappa_h \left( \|\pp_x T\|_2^2
+\|\pp_y T\|_2^2\right) + \kappa_v \|\pp_z T\|_2^2 \right]  \; |T|^4\;
dxdydz  +\aa \kappa_v \|T(z=0)\|_6^6 \\
&&\hskip-.165in = \int_{\Om} Q |T|^4 T \; dxdydz  -\int_{\Om}
\left( uT_x+vT_y+wT_z \right) |T|^4 T \; dxdydz.
\end{eqnarray*}
By integration by parts,  and using (\ref{SA-4})
and the boundary conditions (\ref{B-1})-(\ref{B-3})
we get
\begin{eqnarray}
&&\hskip-.065in   \int_{\Om}
\left( uT_x+vT_y+wT_z \right) |T|^4 T \; dxdydz =0.
\label{DT-10}
\end{eqnarray}
As a result of the above we conclude
\begin{eqnarray*}
&&\hskip-.168in \frac{1}{6} \frac{d \| T \|_{6}^{6} }{d t} +
 5 \int_{\Om} \left[ \kappa_h \left( \|\pp_x T\|_2^2
+\|\pp_y T\|_2^2\right) + \kappa_v \|\pp_z T\|_2^2 \right]  \; |T|^4\;
dxdydz  +\aa \kappa_v \|T(z=0)\|_6^6 \\
&&\hskip-.165in = \int_{\Om} Q |T|^4 T \; dxdydz \leq \|Q\|_2
\|T\|_{10}^5 \leq C \|Q\|_2 \left( \|T\|_6^2 \; \|\nabla T^3\|+
\|T\|_{6}^5 \right).
\end{eqnarray*}
By the Cauchy--Schwarz inequality we get
\begin{eqnarray*}
&&\hskip-.168in
\frac{d \| T \|_{6}^{6} }{d t} +
  \int_{\Om} \left[ \kappa_h \left( \|\pp_x T\|_2^2
+\|\pp_y T\|_2^2\right) + \kappa_v \|\pp_z T\|_2^2 \right]  \; |T|^4\;
dxdydz  +\aa \kappa_v \|T(z=0)\|_6^6 \\
&&\hskip-.165in = \int_{\Om} Q |T|^4 T \; dxdydz \leq  C \|Q\|_2^2  \|T\|_6^4 + \|Q\|_2
\|T\|_{6}^5  \leq  C \|Q\|_2^2  \|T\|_6^4 +
\|T\|_{6}^6.
\end{eqnarray*}
Thus, from the above and (\ref{P-6}), we have
\begin{eqnarray*}
&&\hskip-.168in
\frac{d \| T \|_{6}^{2} }{d t}  \leq  C \|Q\|_2^2  +
\|T\|_{6}^2  \leq  C  \left[  \|Q\|_2^2  +\|T\|_2^2+ \kappa_h \left( \|\pp_x T \|_2^2
+\|\pp_y T \|_2^2\right) + \kappa_v \|\pp_z T \|_2^2 \right].
\end{eqnarray*}
By integrating the above inequality and using (\ref{T-2}) and (\ref{T-2I}),  we get
\begin{eqnarray}
&&\hskip-.68in \| T (t)\|_6^2  \leq C \left[ (1+\|Q\|_{2}^2) \;(1+ t) + \|T_0\|_{H^1(\Om)}^2 \right].     \label{K-T}
\end{eqnarray}

\subsection{$H^1$ estimates}
Recall again that  $T \in L_{\mbox
{loc}}^{\infty}([0, \mathcal{T}_*), H^1(\Om)) \cap
L_{\mbox {loc}}^2([0, \mathcal{T}_*), H^2(\Om))$,
and since, by Theorem \ref{T-WEAK}, equation
(\ref{SA-5}) holds in $L_{\mbox
{loc}}^2([0,\mathcal{T}_*);L^2(\Om))$ we can take
the  inner product of the equation  (\ref{SA-5})
with $-\kappa_h\left( T_{xx} + T_{yy}\right)
 - \kappa_v T_{zz}$, in
$L^2(\Om)$,  and use (\ref{DER}) to obtain, thanks
to  a Lemma of Lions-Magenes concerning the
derivative of functions with values in Banach space
(cf. Chap. III-p.169- \cite{TT84}),
\begin{eqnarray*}
&&\hskip-.268in \frac{1}{2} \frac{d}{dt} \left[\kappa_h
\left( \|\pp_x T\|_2^2 +\|\pp_y T\|_2^2\right)
+ \kappa_v \left( \|\pp_z T\|_2^2 +\aa \|T(z=0)\|_2^2 \right)\right]
+ \|\kappa_h\left(  T_{xx} +  T_{yy}\right)
 + \kappa_v T_{zz}\|_2^2
  \\
&&\hskip-.265in = - \int_{\Om} Q \left[ \kappa_h\left(  T_{xx} +  T_{yy}\right)
+ \kappa_v T_{zz}\right] \; dxdydz  + \int_{\Om} \left( u
\pp_x T +v \pp_y T+ w \pp_z T \right) \left[ \kappa_h\left(  T_{xx} +  T_{yy}\right)
 + \kappa_v T_{zz}\right] \; dxdydz \\
&&\hskip-.268in
\leq  \left[ \|Q\|_2 + \left(\|u\|_6+\|v\|_6 +\|w\|_6\right) \|\nabla T\|_3 \right] \left\| \kappa_h\left(  T_{xx} +  T_{yy}\right)
 + \kappa_v T_{zz}\right\|_2 \\
&&\hskip-.268in
\leq  \left[ \|Q\|_2 + C \|T\|_6^{3/2} \left\| \kappa_h\left(  T_{xx} +  T_{yy}\right)
 + \kappa_v T_{zz}\right\|_2^{1/2} \right] \left\| \kappa_h\left(  T_{xx} +  T_{yy}\right)
 + \kappa_v T_{zz}\right\|_2.
\end{eqnarray*}
By the Cauchy--Schwarz and Young's inequalities we
obtain
\begin{eqnarray*}
&&\hskip-.68in  \frac{d}{dt} \left[\kappa_h
\left( \|\pp_x T\|_2^2 +\|\pp_y T\|_2^2\right)
+ \kappa_v \left( \|\pp_z T\|_2^2 +\aa \|T(z=0)\|_2^2 \right)\right]   + \|\kappa_h\left(  T_{xx} +  T_{yy}\right)
 + \kappa_v T_{zz}\|_2^2
  \\
&&\hskip-.68in
\leq  C \|Q\|_2^2 + C \|T\|_6^6.
\end{eqnarray*}
By Gronwall,  we get
\begin{eqnarray}
&&\hskip-.68in
\kappa_h
\left( \|\pp_x T (t)\|_2^2 +\|\pp_y T(t)\|_2^2\right)
+ \kappa_v \left(\|\pp_z T(t)\|_2^2 +\aa \|T(z=0)(t)\|_2^2 \right) \nonumber \\
&&\hskip-.58in
+\int_0^t \|\kappa_h\left(  T_{xx}(s) +  T_{yy}(s)\right)
 + \kappa_v T_{zz}(s)\|_2^2 \; ds   \nonumber \\
&&\hskip-.68in
\leq C (1+\|Q\|_{2}^2 +\|T\|_6^6 ) \; t + \|
T_0\|_{H^1(\Om)}    \nonumber   \\
&&\hskip-.68in
 \leq C (1+\|Q\|_{2}^2 ) \; t + C \left[ (1+\|Q\|_{2}^2) \;(1+ t) + \|T_0\|_{H^1(\Om)}^2 \right]^3 \; t+ \|
T_0\|_{H^1(\Om)} =: K_v(t).     \label{K-H1}
\end{eqnarray}
Thus,
\[
\limsup_{t \to \mathcal{T}_*^{-}} \|  T \|_{H^1(\Om)} = K_v(\mathcal{T}_*).
\]
This contradicts the assumption that $\mathcal{T}_*$ is finite,  therefore,
$\mathcal{T}_*= \infty$, and the solution $(u,v,w,p, T)$ exists globally in time.

\subsection{Uniqueness of the strong solution and continuous dependence on initial data}

Next,  we show the continuous dependence on the initial data and
the  the uniqueness of the strong solutions. Let $(u_1, v_1, w_1, p_1, T_1)$ and
$(u_2, v_2, w_2, p_2, T_2)$ be two strong solutions of the system
(\ref{SA-1})--(\ref{SA-5}) with corresponding  initial data $(T_0)_1$ and
$(T_0)_2$, respectively. Denote by $u=u_1-u_2, v=v_1-v_2, w=w_1-w_2, p
=p_1 -p_2$ and $ \tt = T_1-T_2.$ It is clear that
\begin{eqnarray}
&&\hskip-.8in
\ee \, u - f \, v + p_x   = 0,    \label{USA-1}  \\
&&\hskip-.8in
 \ee\, v + f \, u  + p_y  = 0,   \label{USA-2}   \\
&&\hskip-.8in
\dd \, w + p_z    = \tt , \label{USA-3}  \\
&&\hskip-.8in
u_x +v_y + w_z =0   \label{USA-4} \\
&&\hskip-.8in \pp_t \tt -\kappa_h\left(  \tt_{xx} +  \tt_{yy}\right)
 - \kappa_v \tt_{zz} + u_1 \tt_x + v_1 \tt_y + w_1 \tt_z
 + u \pp_x T_2 + v\pp_y T_2 + w \pp_z T_2 = 0 \,, \label{USA-5}
\end{eqnarray}
and $(u,v,w)$ and $\tt$ satisfy boundary conditions
(\ref{B-1})--(\ref{B-3}). By Theorem \ref{T-WEAK}
and  Theorem \ref{T-MAIN} equation (\ref{USA-5})
holds in $L^2([0, \mathcal{T}];L^2(\Om))$ and $\tt
\in L^{\infty}([0, \mathcal{T}), H^1(\Om)) \cap
L^2([0, \mathcal{T}), H^2(\Om))$, for all
$\mathcal{T}>0$. Therefore,  by taking the inner
product of equation (\ref{USA-5}) with $\tt$ in
$L^2(\Om)$, and using boundary conditions
(\ref{B-1})--(\ref{B-3}), we get
\begin{eqnarray*}
&&\hskip-.268in \frac{1}{2} \frac{d \|\tt\|_2^2}{dt}
+
 \kappa_h
\left( \|\pp_x \tt \|_2^2 +\|\pp_y \tt \|_2^2\right)
+ \kappa_v \|\pp_z \tt \|_2^2 +\aa \|\tt(z=0)\|_2^2  \\
&&\hskip-.265in =   - \int_{\Om} \left[ u_1 \tt_x + v_1 \tt_y + w_1 \tt_z
 + u (T_2)_x + v(T_2)_y + w (T_2)_z
\right]\; \tt  \; dxdydz.
\end{eqnarray*}
By integration by parts and again boundary conditions (\ref{B-1})--(\ref{B-3}), we get
\begin{eqnarray}
&&\hskip-.065in   - \int_{\Om} \left[ u_1 \tt_x + v_1 \tt_y + w_1 \tt_z
\right]\; \tt  \; dxdydz =0.
\label{DUT-1}
\end{eqnarray}
Notice that
\begin{eqnarray*}
&&\hskip-.68in \left| \int_{\Om} \left[ u (T_2)_x + v(T_2)_y + w (T_2)_z
\right]\; \tt  \; dxdydz \right| \leq C \|\nabla T_2\|_2 \left(\|u\|_4+\|v\|_4+\|w\|_4 \right)\, \|\tt\|_4     \\
&&\hskip-.68in
 \leq C \|\nabla T_2\|_2 \left(\|u\|_2^{1/4} \|u\|_{H^1}^{3/4}+\|v\|_2^{1/4} \|v\|_{H^1}^{3/4}+\|w\|_2^{1/4} \|w\|_{H^1}^{3/4} \right)\,
 \|\tt\|_2^{1/4} \|\tt\|_{H^1}^{3/4}    \\
&&\hskip-.68in
 \leq C \|\nabla T_2\|_2
 \|\tt\|_2^{1/2} \|\tt\|_{H^1}^{3/2}
\leq C
\|\nabla T_2\|_2 \left( \|\tt\|_2^{2} + \|\tt\|_2^{1/2}   \|\nabla \tt\|_2^{3/2} \right).
\end{eqnarray*}
Thus,
\begin{eqnarray*}
&&\hskip-.268in
\frac{1}{2} \frac{d \|\tt\|_2^2}{dt}
+
 \kappa_h
\left( \|\pp_x \tt \|_2^2 +\|\pp_y \tt \|_2^2\right)
+ \kappa_v \|\pp_z \tt \|_2^2 +\aa \|\tt(z=0)\|_2^2  \\
&&\hskip-.265in \leq C
\|\nabla T_2\|_2 \left( \|\tt\|_2^{2} + \|\tt\|_2^{1/2}   \|\nabla \tt\|_2^{3/2} \right).
\end{eqnarray*}
By Young's inequality, we get
\begin{eqnarray*}
&&\hskip-.268in
 \frac{d \|\tt\|_2^2}{dt}
+
 \kappa_h
\left( \|\pp_x \tt \|_2^2 +\|\pp_y \tt \|_2^2\right)
+ \kappa_v \|\pp_z \tt \|_2^2 +\aa \|\tt(z=0)\|_2^2  \\
&&\hskip-.265in \leq C \|\nabla T_2\|_2^4 \|\tt\|_2^{2}.
\end{eqnarray*}
Thanks to Gronwall inequality, we obtain
\begin{eqnarray*}
&&\hskip-.68in \|\tt(t)\|_2^2
\leq\|\tt(t=0)\|_2^2  e^{ C  \int_0^t  \|\nabla T_2
(s)\|_2^4\; ds}.
\end{eqnarray*}
Since $T_2$ is a strong solution,  we have by
virtue of (\ref{K-H1})
\begin{eqnarray}
&&\hskip-.68in \|\tt(t)\|_2^2 \leq \|\tt(t=0)\|_2^2\,
e^{ C \int_0^tK_v^2(s) ds},    \label{CCC}
\end{eqnarray}
where the value of $T_0$ in the definition of $K_v$
in (\ref{K-H1}) is replaced by $T_2(0)$. As a
result, the above inequality proves the continuous
dependence of the solutions on the initial data. In
particular, when $\tt(t=0)=0$, we have $\tt(t)=0,$
and consequently also $u(t)=v(t)=w(t)=0$, for all
$t\ge 0$. Therefore, the strong solution is unique.

\end{proof}

\noindent
\section*{Acknowledgements}
 This work was
supported in part by the NSF grants
no.~DMS-0709228, no.~DMS-0708832, and
no.~DMS-1009950, and by the Alexander von Humboldt
Stiftung/Foundation and  the Minerva
Stiftung/Foundation  (E.S.T.). The authors are
thankful to the kind hospitality of  the Institute
for Mathematics and its Applications (IMA),
University of Minnesota, and E.S.T. also thankful
to the warm hospitality of the Freie
Universit\"{a}t--Berlin,  where part of this work
was completed.

\end{document}